\documentclass[a4paper,11pt,english]{amsart}
\usepackage[applemac]{inputenc}
\usepackage[english]{babel}
\usepackage{amsfonts}
\usepackage{amsmath} 
\usepackage{amsthm}
\usepackage{hyperref} 
\usepackage{latexsym}
\usepackage{mathrsfs}
\usepackage{array}
\usepackage{amssymb}
\usepackage{enumerate}
\usepackage{graphicx}
\usepackage{color}
\usepackage{stmaryrd}
\newcommand{\R}{\mathbb  R}
\newcommand{\C}{\mathbb  C}

\renewcommand{\epsilon}{\varepsilon}

\newcommand{\Z}{  \mathbb{Z}   }

\newcommand{\dis}{\displaystyle}

\newcommand{\ov}{  \overline  }
\renewcommand{\a}{  \alpha   }

\renewcommand{\phi}{  \varphi  }

\renewcommand{\>}{  \rangle   }

\numberwithin{equation}{section}
\theoremstyle{plain}

\newtheorem{prop}{Proposition}[section]

\def\beq{\begin{equation}}   \def\eeq{\end{equation}}
\def\bea{\begin{eqnarray}}  \def\eea{\end{eqnarray}}

\renewcommand{\theequation}{\thesection.\arabic{equation}}
\newcounter{hran} \renewcommand{\thehran}{\thesection.\arabic{hran}}

\def\bmini{\setcounter{hran}{\value{equation}}
    \refstepcounter{hran}\setcounter{equation}{0}
    \renewcommand{\theequation}{\thehran\alph{equation}}\begin{eqnarray}}

\def\bminiG#1{\setcounter{hran}{\value{equation}}
\refstepcounter{hran}\setcounter{equation}{-1}
\renewcommand{\theequation}{\thehran\alph{equation}}
\refstepcounter{equation}\label{#1}\begin{eqnarray}}

\author{Pierre Germain}
\address{Courant Institute of Mathematical Sciences, 251 Mercer Street, New York 10012-1185 NY, USA
}
\email{pgermain@cims.nyu.edu}
\author{Laurent Thomann}
\address{
Institut  \'Elie Cartan, Universit\'e de Lorraine, B.P. 70239,
F-54506 Vandoeuvre-l\`es-Nancy Cedex, France}

\email{laurent.thomann@univ-lorraine.fr}

\title[High frequency limit of the LLL equation]{On the high frequency limit of the LLL equation}
 
\begin{document}

\begin{abstract}
We derive heuristically an integro-differential equation, as well as a shell model, governing the dynamics of the Lowest Landau Level equation in a high frequency regime. 
\end{abstract}

\thanks{P. G. is partially supported by NSF grant DMS-1101269, a start-up grant from the Courant Institute, and a Sloan fellowship.}
\thanks{L.T. is partially supported   by the  grant  ``ANA\'E'' ANR-13-BS01-0010-03.}
\maketitle

\section{Introduction}

\subsection{The LLL equation}

The   Lowest Landau Level evolution equation (LLL for short) reads
\begin{equation}
\label{LLLu}
i \partial_t u = \Pi\big( |u|^2 u\big),
\end{equation}
where the unknown function $u = u(t,z)$ is a map from $\mathbb{R} \times \mathbb{C}$ to $\mathbb{C}$, and  where $\Pi$ is the orthogonal projector in $L^2(\mathbb{C})$ on the Bargmann-Fock space $\mathcal{E} = \{ e^{-\frac{|z|^2}{2}} f(z), \; \mbox{where $f$ is entire}\} \cap L^2 (\mathbb{C})$. Denote by $L$   the Lebesgue measure on $\C$, then the  kernel of $\Pi$ is explicitly given by 
$$
[\Pi u](z) = \frac{1}{\pi} e^{-\frac{|z|^2}{2}} \int_\mathbb{C} e^{\ov  w z - \frac{|w|^2}{2}} u(w) \,dL(w).
$$
The equation \eqref{LLLu} derives from the Hamiltonian
$$
\mathcal{H}(u) = \int_{\mathbb{R}^2} |u(z)|^4 \,L(dz),
$$
defined on $\mathcal{E}$ equipped the symplectic form $\dis \omega(f,g) = \int f \overline{g} \, L(dz)$ defined on $L^2 ( \mathbb{R}^2)$.\bigskip

This equation appears in the context of rapidly rotating Bose-Einstein condensates, and more precisely in the study of vibration modes of condensates, the so called Tkachenko modes, see \cite{BWCB, SCEMC, Son}. 
From a more mathematical  point of view, the variational problem was examined in~\cite{ABN,AS}, and the dynamical properties of~\eqref{LLLu} in\;\cite{Nier,GHT1, GHT2, GGT} for more on the analysis of this equation.\bigskip

The equation \eqref{LLLu} is globally well-posed in $\mathcal{E}$ (see  \cite{Nier}). It is then  natural to study the long-time dynamics, and in particular to decide if  large amounts of  energy can be transfered from low to high frequencies. Such a phenomenon would also have counterparts for the 3-dimensional Gross-Pitaevskii equation (see \cite{HT} for more details).  \bigskip

In this paper we derive, heuristically,  an integro-differential equation which may govern the dynamics of the LLL equation. In the same spirit, we also derive a shell model. Understanding the dynamics of these models (in particular the possible  transfer of energy) may help to catch key features of the  dynamics of the full model.

\subsection{Our Ansatz} As discussed in~\cite{GHT1,GHT2}, a basis of the Bargmann-Fock space which is particularly well suited to the analysis of LLL is provided by
\begin{equation*}
\phi_n(z) = \frac{1}{\sqrt{\pi n!}} z^n e^{-\frac{|z|^2}{2}}, \quad \mbox{with $n \in \mathbb{N}$},
\end{equation*}
so that a natural Ansatz for the high frequency limit  of~\eqref{LLLu} is given by
\begin{equation}\label{defvl}
u = a \phi_0 + v_\lambda \quad \mbox{with} \quad v_\lambda(t,z) = \frac{1}{\sqrt{\lambda}} \sum_{n \in \mathbb{N}} g\left(t, \frac{n}{\lambda} \right) \phi_{n}(z) = \frac{1}{\sqrt{\lambda}} \sum_{k \in \mathbb{N}/\lambda}g(t,k) \phi_{\lambda k}(z),
\end{equation}
where $g$ is such that $g(\cdot,0) = 0$. The high frequency limit corresponds to $\lambda \to \infty$.

In this limit, it is natural to view $g$ as a continuous function of its second variable; our aim is to derive an effective equation for it.

\subsection{The limiting equation} 

It turns out to be simpler to pass to the limit in the Hamiltonian. 

With the above Ansatz, we will show \underline{heuristically} in Section~\ref{sectionderivation} that, up to a multiplicative constant,
$$
\mathcal{H}(u) \overset{\lambda \to \infty}{\sim}  \mathfrak{h}_\lambda(a,g): = |a|^4  + 4 \int_0^{+\infty} \Big\vert   \frac{ag(s)}{2^{\lambda s/2}}   +  K \left( \frac{s}{\lambda}\right)^{\frac{1}{4}}     {g^2 \left( \frac{s}{2} \right)} \Big\vert^2 \,ds,
$$
with $K =2^{-3/4} \pi^{1/4}$. As $\lambda \to \infty$, this leads to the dynamics $i \dot a  = \frac{\partial \mathfrak{h}_\lambda}{\partial \overline a}$ and $i \dot g = \frac{1}{\lambda} \frac{\partial \mathfrak{h}_\lambda}{\partial \overline g}$, where the $\frac{1}{\lambda}$ factor is due to our Ansatz. In other words,
$$
\left\{
\begin{array}{l}
\displaystyle i \dot{a} = 2 |a|^2 a + 4  a \int_0^{+\infty}  \frac{|g(s)|^2}{2^{\lambda s}} \,ds + 4K \int_0^{+\infty}  \frac{\overline{g(s)}}{2^{\lambda s/2}}\left( \frac{s}{\lambda}\right)^{\frac{1}{4}} g^2 \left( \frac{s}{2} \right)\,ds \\
\displaystyle i \dot{g}(s) = \frac{1}{\lambda} \left[ 4 |a|^2 \frac{g(s)}{2^{\lambda s}} + 4K \left( \frac{s}{\lambda}\right)^{\frac{1}{4}} \frac{\overline{a}}{2^{\lambda s/2}} g^2 \left( \frac{s}{2} \right) + 8K\left( \frac{2s}{\lambda}\right)^{\frac{1}{4}} \frac{a}{2^{\lambda s}} \overline{g(s)} g(2s) + 16K^2  \left( \frac{2s}{\lambda}\right)^{\frac{1}{2}} |g(s)|^2 g(s) \right].
\end{array}
\right. 
$$
One of the very attractive peculiarities of this system is that a given frequency $s$ only interacts with the frequencies $\frac{s}{2}$ or $2s$. In particular, it seems natural, in order to further simplify the system, to assume that $g$ is supported on $2^{\mathbb{Z}}$. With this aim in mind, observe that the Ansatz
$$
g = \sum_{j \in \mathbb{Z}} g_j \mathbf{1}_{[2^j,2^j(1+\epsilon)]}
$$
leads, as $\epsilon \to 0$, to dynamics given by a shell model:
$$
\left\{
\begin{array}{l}
\displaystyle i \dot{a} = 2 |a|^2 a + 4  a \epsilon \sum_{j\in \Z} 2^j \frac{|g_j|^2}{2^{\lambda{2^j}}} + 4K \epsilon  \sum_{j\in \Z} 2^j \frac{\overline{g_j}}{2^{\lambda 2^{j-1}}}\left( \frac{2^j}{\lambda}\right)^{\frac{1}{4}} g_{j-1}^2 \\
\displaystyle i \dot{g_j} =  \frac{1}{\lambda} \left[ 4  |a|^2 \frac{g_j}{2^{\lambda 2^j}} + 4K \left( \frac{2^j}{\lambda}\right)^{\frac{1}{4}} \frac{\overline{a}}{2^{\lambda 2^{j-1}}} g_{j-1}^2 + 8K\left( \frac{2^{j+1}}{\lambda}\right)^{\frac{1}{4}} \frac{a}{2^{\lambda 2^j}} \overline{g_j} g_{j+1} + 16K^2  \left( \frac{2^{j+1}}{\lambda}\right)^{\frac{1}{2}} |g_j|^2 g_j \right].
\end{array}
\right.
$$
As a shell model, the above equation only retains information about dyadic frequency ranges. This particularly simple structure makes its simulation computationally much more accessible, while conserving essential nonlinear features. 

\subsection{Comparison to other models} \underline{\textit{Shell models}} have been influential in turbulence theory~\cite{G, OY,B}. They also rely on a dyadic decomposition in frequency, except that now $u_j \in \mathbb{C}$  accounts for the whole frequency range $[2^j,2^{j+1}]$. For the popular GOY model, the dynamics is given by
$$
\dot{u}_j = -\nu 2^{2j} u_j + i \left( 2^j \overline{u_{j+2}} \overline{u_j} - b 2^{j-1} \overline{u_{j+1}} \overline{u_{j-1}} - (1-b) \overline{u_{j-1}} \overline{u_{j-2}} \right)
$$
(where $\nu>0$ is the viscosity, and $b$ a parameter depending on the system). This system of ODE enjoys two conserved quantities, and has been successful in reproducing some aspects of hydrodynamic turbulence. However, its derivation is only phenomenological, and it lacks a Hamiltonian structure for $\nu = 0$.

The model that we propose eschews these two shortcomings, and we hope it can lead to some insight into weak turbulence.

\bigskip

\noindent
\underline{\textit{The toy model}} proposed in~\cite{CKSTT}, see also~\cite{CMOS,HM}, derives from the Hamiltonian $\sum_{j=1}^\infty \Big(\frac{1}{4} |b_j|^4 - \mathfrak{Re} \overline{b_j^2} b_{j-1}^2\Big)$. It reads
$$
i \dot{b}_j = |b_j|^2 b_j + 2 b_{j-1}^2 \overline{b_j} + 2 b_{j+1}^2 \overline{b_j}.
$$
It shares many features with the model we propose here, both as far as the structure of the system goes; and in the aim of modeling weakly turbulent flow. But it corresponds to a very different context: it can be seen as a subsystem of NLS on the 2-torus, in the case when the data is sparse in Fourier, and satisfies a certain condition of "separation between generations".

\bigskip

\noindent
\underline{\textit{A model for growth of high frequencies in the weakly turbulent regime}} has been derived in the work~\cite[Section\;5.1]{EV}. More precisely, this reference deals with the kinetic wave equation (or Kolmogorov-Zakharov equation) and the idea pursued is to model the interaction of a (Bose-Einstein) condensate at frequency zero with a high frequency perturbation. This is of course similar to the model proposed here, where part of the solution is in the ground state ($a \phi_0$) while the other one is at high frequency\;($v_\lambda$). However, since it is derived from the kinetic wave equation, it is an equation of kinetic type instead of a Hamiltonian system.

\section{Derivation of the asymptotic Hamiltonian}

\label{sectionderivation}

\subsection{A few more words on LLL} 
Recall that
$$
\phi_n(z) = \frac{1}{\sqrt{\pi n!}} z^n e^{-\frac{|z|^2}{2}}, \quad \mbox{with $n \in \mathbb{N}$}
$$
gives an orthonormal basis of the Bargmann-Fock space (see~\cite{Zhu}). Expanding $u$ in this basis
$$
u(t,z) = \sum_{n=0}^\infty c_n(t) \phi_n(z),
$$
and accelerating time by a factor $\frac{\pi}{2}$, a simple
computation~\cite{GHT1} gives that the equation~\eqref{LLLu} becomes
$$
i \dot{c_n} = 2 \sum_{\substack{k+\ell=m+n \\ k,\ell,m,n \in \mathbb{N}}} \frac{(k+\ell)!}{2^{k+\ell} \sqrt{k!\, \ell!\, m!\, n!}} c_k c_{\ell} \overline{c_m}.
$$
or
$$
i \dot{c_n} = \frac{\partial \mathcal{H}(c)}{\partial \overline{c_n}}
$$
after defining
$$
\mathcal{H}(e,f,g,h) = \sum_{\substack{k+\ell=m+n \\ k,\ell,m,n \in \mathbb{N}}} \frac{(k+\ell)!}{2^{k+\ell} \sqrt{k!\, \ell!\, m!\, n!}} e_k f_{\ell} \overline{g_m} \overline{h_n}, \quad \mbox{and} \quad \mathcal{H}(c) = \mathcal{H}(c,c,c,c).
$$

\subsection{Expanding the Hamitonian}\label{sect22} Recall the definition of\;$v_{\lambda}$ in \eqref{defvl}. Then we write
\begin{equation*}
\begin{split}
 \mathcal{H}(a\phi_0 + v_\lambda) = &\underbrace{ |a|^4 \mathcal{H}(\phi_0)}_{\mathcal{H}_0^\lambda(a,a,a,a)} + \underbrace{2 \left[ \mathcal{H}(v_\lambda,a\phi_0,a\phi_0,a\phi_0) +  \mathcal{H}(a\phi_0,a\phi_0,v_\lambda,a\phi_0) \right]}_{\mathcal{H}_1^\lambda(a,a,a,g)}\\
& \quad + \underbrace{\left[ \mathcal{H}(v_\lambda,v_\lambda,a\phi_0,a\phi_0) + \mathcal{H}(a\phi_0,a\phi_0,v_\lambda,v_\lambda) + 4 \mathcal{H}(v_\lambda,a\phi_0,v_\lambda,a\phi_0) \right]}_{\mathcal{H}_2^\lambda(a,a,g,g)} \\
& \quad + \underbrace{ 2\left[ \mathcal{H}(v_\lambda,v_\lambda,v_\lambda,a\phi_0) +   \mathcal{H}(a\phi_0,v_\lambda,v_\lambda,v_\lambda) \right]}_{\mathcal{H}_3^\lambda(a,g,g,g)} + \underbrace{\mathcal{H}(v_\lambda)}_{\mathcal{H}_4^\lambda(g,g,g,g)}.
\end{split}
\end{equation*}

\subsection{Cancellation of  \texorpdfstring{$\mathcal{H}^\lambda_1(a,a,a,g)$}{H1}} Since $g(0)=0$, one sees immediately that $\mathcal{H}^{\lambda}_1(a,a,a,g)=0$.

\subsection{Equivalent for \texorpdfstring{$\mathcal{H}^{\lambda}_2(a,a,g,g)$}{H2}} 
It is easy to see that
$$
\mathcal{H}_2^\lambda(a,a,g,g) = \frac{2}{\lambda} \mathfrak{Re} \;\overline{a}^2\underbrace{ g^2(0)}_0 + 4 |a|^2 \sum_{k \in \mathbb{N}/\lambda} \frac{(\lambda k)!}{2^{\lambda k} (\lambda k)! \lambda} |g(k)|^2.
$$
Thus, as $\lambda \to \infty$, taking the continuous limit gives
$$
\mathcal{H}_2^\lambda(a,a,g,g) \sim 4 |a|^2 \int_{0}^{+\infty} \frac{1}{2^{\lambda x}} |g(x)|^2 \,dx.
$$

\subsection{Equivalent for $\mathcal{H}_3^\lambda(a,g,g,g)$} By the formula giving $\mathcal{H}$,
$$
\mathcal{H}_3^\lambda(a,g,g,g) = 4 \,\mathfrak{Re}\, a \sum_{\substack{k = m+n \\ k,m,n \in \mathbb{N}/\lambda}} \frac{(\lambda k)!}{2^{\lambda k} \sqrt{(\lambda k)!(\lambda m)!(\lambda n)!}} \frac{1}{\lambda^{3/2}} g(k) \overline{g(m) g(n)}.
$$
By Stirling's formula
$$
\mathcal{H}_3^\lambda(a,g,g,g) \sim  \frac{2^{7/4} }{\pi^{1/4}\lambda^{7/4}}\,\mathfrak{Re}\, a \sum_{k = m+n}  \left( \frac{k}{mn} \right)^{1/4} \left( \frac{k^{k}}
{2^{2k} m^{m} n^{n}} \right)^{\lambda/2} g(k) \overline{g(m) g(n)}.
$$
Take now the continuous limit ($k$ and $m$ becomes the continuous variables $x$ and $y$ respectively):
$$
\mathcal{H}_3^\lambda \sim \frac{2^{7/4}  \lambda^{1/4}}{\pi^{1/4}} \,\mathfrak{Re}\, a \int_{x\in \R_+} \int_{0<y<x} \left( \frac{x}{y(x-y)} \right)^{1/4} \left( \frac{x^x}{2^{2x} y^y (x-y)^{x-y}} \right)^{\lambda/2}  g(x) \overline{g(y) g(x-y)} \,dy\,dx.
$$
Changing the integration variable to $\theta = \frac{y}{x}$, the above becomes
\begin{align*}
\mathcal{H}_3^\lambda & \sim  \frac{2^{7/4}  \lambda^{1/4}}{\pi^{1/4}}  \, \mathfrak{Re} \, a \int_{x>0} \int_{0<\theta<1} \left( \frac{1}{x\theta(1-\theta)} \right)^{1/4} \left( \frac{1}{4 \theta^\theta (1-\theta)^{1-\theta}} \right)^{\lambda x/2} g(x) \overline{g(\theta x) g((1-\theta)x)} \,x\,d\theta\,dx \\
& =  \frac{2^{7/4}  \lambda^{1/4}}{\pi^{1/4}} \, \mathfrak{Re}\, a \int_{x>0} \int_{0<\theta<1} \left( \frac{x^3}{\theta(1-\theta)} \right)^{1/4} \left( \frac{1}{2} \right)^{\lambda x } \psi(\theta)^{\lambda x/2} g(x) \overline{g(\theta x) g((1-\theta)x)} \,d\theta\,dx,
\end{align*}
where we denote
$$
\psi(\theta) = \frac{1}{\theta^\theta (1-\theta)^{1-\theta}}.
$$
Let us take a moment to study this function, which is defined on $(0,1)$. It is maximal and equal to\;2 for $\theta = 1/2$, and approaches 1 as $\theta$ approaches 0 or 1. Furthermore, its second derivative at $1/2$ equals\;$8$, and thus $\psi(1/2+\delta) = 2 -4\delta^2 + O(\delta^3)$. As $\lambda \to \infty$, it is clear that, for a smooth function $F$, a positive real $\alpha$, and $\lambda \to \infty$,
\begin{equation}
\label{psilambda}
\int_{0<\theta<1} \psi(\theta)^{\alpha \lambda} F(\theta) \,d\theta \sim F\left(\frac{1}{2}\right) \int_{0}^{1} 2^{\alpha \lambda} e^{- 2 \alpha \lambda (\theta - 1/2)^2}\,d\theta \sim F\left(\frac{1}{2}\right) \frac{2^{\alpha \lambda} \sqrt{\pi}}{\sqrt{2 \alpha \lambda}}.
\end{equation}
Using this in the equivalent for $\mathcal{H}_3$ gives
$$
\mathcal{H}_3^\lambda(a,g,g,g) \sim 2^{9/4} \pi^{1/4} \lambda^{-1/4} \, \mathfrak{Re} \,a \int_{x>0} x^{1/4} \left( \frac{1}{2} \right)^{\lambda x/2} \overline{ g \left(\frac{x}{2} \right)}^2 g(x) \,dx.
$$

\subsection{Equivalent for  \texorpdfstring{$\mathcal{H}$}{H}}
By the formula giving $\mathcal{H}$,
$$
\mathcal{H}(v_\lambda) = \sum_{\substack{k+\ell= m+n \\ k,\ell,m,n \in \mathbb{N}/\lambda}}  \frac{\big(\lambda (k+\ell)\big)!}{2^{\lambda (k+\ell)} \sqrt{(\lambda k)!(\lambda \ell)!(\lambda m)!(\lambda n)!}} \frac{1}{\lambda^{2}} g(k) g(\ell) \overline{g(m) g(n)}.
$$
By Stirling's formula
\begin{align*}
\mathcal{H}(v_\lambda) & \sim  \frac{1}{\sqrt{2\pi}} \sum_{k+\ell = m+n}\frac{1}{\lambda^{5/2}} \left( \frac{(k+\ell)^2}{k\ell mn} \right)^{1/4} \left( \frac{(k+\ell)^{k+\ell}}
{2^{k+\ell} \sqrt{k^{k}\ell^{\ell}m^{m} n^{n}}} \right)^\lambda g(k) g(\ell) \overline{g(m) g(n)} \\
&=  \frac{1}{\sqrt{2 \pi}} \frac{1}{\lambda^{5/2}} \sum_{\substack{S \\ 0\leq k \leq S \\ 0\leq m \leq S}} \left( \frac{S^2}{k(S-k)m(S-m)} \right)^{1/4} \left[ \frac{1}{2} \psi \left( \frac{k}{S} \right) \right]^{\lambda S/2} \left[ \frac{1}{2} \psi \left( \frac{m}{S} \right)\right]^{\lambda S/2} \\
& \qquad \qquad\qquad\qquad\qquad\qquad\qquad\qquad\qquad\qquad\qquad\qquad\qquad g(k) g(S-k) \overline{g(m) g(S-m)}
\end{align*}
(where we set $S = k+\ell = m+n$).
Taking the continuous limit (the discrete variables $k,m,S$ becoming $x,y,s$ respectively) and subsequently changing variables to $\theta = x/s$ and $\eta = y/s$ gives
\begin{align*}
\mathcal{H}(v_\lambda) & \sim \frac{\sqrt{\lambda}}{\sqrt{2 \pi}} \int_{s>0} \int_{0<x<s} \int_{0<y<s} \left( \frac{s^2}{x(s-x)y(s-y)} \right)^{1/4} \left[ \frac{1}{2} \psi \left( \frac{x}{s} \right) \right]^{\lambda s/2} \left[\frac{1}{2} \psi \left( \frac{y}{s} \right) \right]^{\lambda s/2} \\
& \qquad \qquad\qquad\qquad\qquad\qquad\qquad\qquad\qquad\qquad\qquad g(x)g(s-x)\overline{g(y)g(s-y)}\,dx\,dy\,ds \\
& = \frac{\sqrt{\lambda}}{\sqrt{2 \pi}}\int_{s>0} \int_{0<\theta<1} \int_{0<\eta<1} \frac{1}{\sqrt{s}} \left( \frac{1}{\theta(1-\theta)\eta(1-\eta)} \right)^{1/4} \left[ \frac{1}{2} \psi \left( \theta \right) \right]^{\lambda s/2} \left[\frac{1}{2} \psi \left(\eta \right) \right]^{\lambda s/2} \\
& \qquad \qquad\qquad\qquad\qquad\qquad\qquad\qquad\qquad\qquad\qquad g(\theta s)g((1-\theta)s)\overline{g(\eta s)g((1-\eta)s)}\,s\,d\theta\,s\,d\eta\,ds.
\end{align*}
By~(\ref{psilambda}), this is further equivalent to
$$
\mathcal{H}(v_\lambda) \sim \frac{\sqrt{2\pi}}{\sqrt{\lambda}} \int \sqrt{s} \left| g \left( \frac{s}{2} \right) \right|^4 \,ds .
$$

\subsection{What is the relevant expansion for $\mathcal{H}$?}

Recall that we expanded in Section\;\ref{sect22}
$$
\mathcal{H}(u) = \mathcal{H}_0^\lambda (a,a,a,a) + \mathcal{H}_1^\lambda(a,a,a,g) + \mathcal{H}_2^\lambda(a,a,g,g) + \mathcal{H}_3^\lambda(a,g,g,g) + \mathcal{H}_4^\lambda(g,g,g,g),
$$
where, for instance, $\mathcal{H}^{\lambda}_1(a,a,a,g)$, is (real) trilinear in $a$ and linear in $g$. We then proved that
\begin{align*}
& \mathcal{H}_0^\lambda(a,a,a,a) = \mathfrak{h}_0^\lambda(a)=|a|^4 \\
& \mathcal{H}_1^\lambda(a,a,a,g) = 0 \\
& \mathcal{H}_2^\lambda(a,a,g,g) \overset{\lambda \to \infty}{\sim} \mathfrak{h}_2^\lambda(a,a,g,g) = 4 |a|^2 \int_{0}^{+\infty} \frac{1}{2^{\lambda x}} |g(x)|^2 \,dx. \\
& \mathcal{H}_3^\lambda(a,g,g,g) \overset{\lambda \to \infty}{\sim} \mathfrak{h}_3^\lambda(a,g,g,g) = 2^{9/4} \pi^{1/4} \lambda^{-1/4} \, \mathfrak{Re} \,a \int_{x>0} x^{1/4} \left( \frac{1}{2} \right)^{\lambda x/2} \overline{ g \left(\frac{x}{2} \right)}^2 g(x) \,dx. \\
& \mathcal{H}^\lambda_4(g,g,g,g) \overset{\lambda \to \infty}{\sim} \mathfrak{h}^\lambda_4 (g,g,g,g) =  \frac{\sqrt{2\pi}}{\sqrt{\lambda}} \int_{0}^{+\infty} \sqrt{s} \left| g \left( \frac{s}{2} \right) \right|^4 \,ds.
\end{align*}
How can we combine these different asymptotics, or in other words: what is the relevant expansion for $\mathcal{H}$?
\begin{itemize}
\item Keeping only terms of size $O(1)$ (as far as their dependence on $\lambda$ goes) leads to dropping $\mathcal{H}_2^\lambda$, $\mathcal{H}_3^\lambda$ and $\mathcal{H}_4^\lambda$, which gives the uninteresting model $\dis \mathcal{H}(u) \sim \frac{\pi}8 |a|^4$.
\item Keeping only terms which are polynomial in $\lambda$ leads to dropping $\mathcal{H}_2^\lambda$ and $\mathcal{H}_3^\lambda$ (which decay exponentially fast in $\lambda$), thus obtaining 
$$
\mathcal{H}(u) \overset{\lambda \to \infty}{\sim} \mathfrak{h}_0^\lambda + \mathfrak{h}_4^\lambda =|a|^4 + \frac{\sqrt{2\pi}}{\sqrt{\lambda}} \int_{0}^{+\infty} \sqrt{s} \left| g \left( \frac{s}{2} \right) \right|^4 \,ds.
$$
This corresponds once again to trivial dynamics, since for all $t\in \R$ we get for the associated dynamical system $\vert a(t)\vert=\vert a(0)\vert$ and $\vert g(t,s)\vert=\vert g(0,s)\vert$.
\item Therefore, to observe a genuine nonlinear interaction between modes, we have to include the contributions of $\mathcal{H}_3^\lambda$ and $\mathcal{H}_4^\lambda$ to get
$$
\mathcal{H}(u) \overset{\lambda \to \infty}{\sim} \mathfrak{h}_0^\lambda + \mathfrak{h}_2^\lambda + \mathfrak{h}_3^\lambda + \mathfrak{h}_4^\lambda = |a|^4  + 4 \int_0^{+\infty} \Big\vert   \frac{ag(s)}{2^{\lambda s/2}}   +  K \left( \frac{s}{\lambda}\right)^{\frac{1}{4}}     {g^2 \left( \frac{s}{2} \right)} \Big\vert^2 \,ds.
$$
Observe that this expansion must be taken with a grain of salt since, for instance, the error in the expansion of $\mathcal{H}^\lambda_4$ is larger than the leading order term in the expansion of $\mathcal{H}^\lambda_2$. It should be understood as follows: we keep the leading contribution to trivial resonances (self interactions at the frequency $s$, given by $\mathfrak{h}^\lambda_4$), as well as the leading contribution to nontrivial resonances (interactions between the frequencies $s$ and $2s$, given by $\mathfrak{h}^\lambda_2$ and $\mathfrak{h}^\lambda_3$).
\end{itemize}

\section{Some remarks on the limiting equation}

We discuss here the dynamics of the evolution equation given by the above Hamiltonian:
\begin{equation}\label{sys}
\left\{
\begin{array}{l}
 i \dot{a} = 2 |a|^2 a + 4  a \int_0^{+\infty}  \frac{|g(s)|^2}{2^{\lambda s}} \,ds + 4K \int_0^{+\infty}  \frac{\overline{g(s)}}{2^{\lambda s/2}}\left( \frac{s}{\lambda}\right)^{\frac{1}{4}} g^2 \left( \frac{s}{2} \right)\,ds, \\
i \dot{g}(s) = \frac{1}{\lambda} \left[ 4 |a|^2 \frac{g(s)}{2^{\lambda s}} + 4K \left( \frac{s}{\lambda}\right)^{\frac{1}{4}} \frac{\overline{a}}{2^{\lambda s/2}} g^2 \left( \frac{s}{2} \right) + 8K\left( \frac{2s}{\lambda}\right)^{\frac{1}{4}} \frac{a}{2^{\lambda s}} \overline{g(s)} g(2s) + 16K^2  \left( \frac{2s}{\lambda}\right)^{\frac{1}{2}} |g(s)|^2 g(s) \right],
\end{array}
\right.
\end{equation}
with $K =2^{-3/4} \pi^{1/4}$, and where $\lambda\geq 1$. Similar conclusions can be drawn about the discrete (in $s$) system. 

\subsection{Local well-posedness}
For $\alpha\geq 0$, define the space $X_\a$ by the norm $\dis \Vert(a,g)\Vert_{\alpha}=\vert a\vert+\sup_{s\geq 0}\<s\>^{\alpha}\vert g(s)\vert $. Then we have the following well-posedness result

\begin{prop}
Let $\a\geq 1/4$ and $(a_0,g_0)\in X_\a$. Then there exists $T>0$ and a unique solution $(a,g)\in \mathcal{C}\big(\big[0,T]; X_\a\big)$ to the system \eqref{sys} with initial conditions $(a_0,g_0)$.
\end{prop}

The local time of existence is $T\leq C  \Vert(a_0,g_0)\Vert^{-2}_{\alpha}$ and the constant $C>0$ is uniform in $\lambda\geq 1$.

\begin{proof}
Define the norm $ \Vert(a,g)\Vert_{T,\alpha}=\dis \sup_{0\leq t\leq T}\Vert(a(t),g(t))\Vert_{\alpha}$. Then,  we consider the space $X=\big\{  (a,g)\in \mathcal{C}\big(\big[0,T]; X_\a\big),\; \Vert(a,g)\Vert_{\alpha}\leq 2   \Vert(a_0,g_0)\Vert_{\alpha}  \big\}$ and the map $\Phi=(\Phi_1,\Phi_2)$ given by
\begin{equation*}
\left\{
\begin{array}{l}
\displaystyle \Phi_1(a,g)(t) = a_0-i\int_0^t\bigg(2i |a|^2 a  +4 a \int_0^{+\infty}  \frac{|g(s)|^2}{2^{\lambda s}} \,ds +4 K \int_0^{+\infty}  \frac{\overline{g(s)}}{2^{\lambda s/2}}\left( \frac{s}{\lambda}\right)^{\frac{1}{4}} g^2 \left( \frac{s}{2} \right)\,ds\bigg)dt', \\
\displaystyle  \Phi_2(a,g)(t) =g_0(s)- \frac{i}{\lambda} \int_0^t\bigg(4 |a|^2 \frac{g(s)}{2^{\lambda s}}  +4K \left( \frac{s}{\lambda}\right)^{\frac{1}{4}} \frac{\overline{a}}{2^{\lambda s/2}} g^2 \left( \frac{s}{2} \right) +\\
\displaystyle \qquad  \qquad \qquad  \qquad  \qquad  \qquad  +8K\left( \frac{2s}{\lambda}\right)^{\frac{1}{4}} \frac{a}{2^{\lambda s}} \overline{g(s)} g(2s) +16K^2  \left( \frac{2s}{\lambda}\right)^{\frac{1}{2}} |g(s)|^2 g(s)\bigg)dt',
\end{array}
\right.
\end{equation*}
 We show that for $T>0$ small enough, the map $\Phi : X \longrightarrow X$ is a contraction, and hence admits a unique fixed point, which is the unique solution to the system\;\eqref{sys}. 
 
 We show that $\Phi$ maps $X$ into itself. For $0\leq t\leq T$, we   get
  \begin{eqnarray*}
  \vert  \Phi_1(a,g)(t)\vert &\leq& \vert a_0\vert+CT \Big( (\sup_{0\leq t\leq T} \vert a(t)\vert)^3+\big(\sup_{0\leq t\leq T} \vert a(t)\vert\big)\big(\sup_{0\leq t\leq T}\sup_{s\geq 0} (\<s\>^{\alpha}\vert g(s)\vert)\big)^2 \int_{0}^{+\infty} \frac{ds}{2^{\lambda s} \<s\>^{2\alpha}}\\
 & &\qquad\qquad\qquad + \big(\sup_{0\leq t\leq T}\sup_{s\geq 0} (\<s\>^{\alpha}\vert g(s))\vert\big)^3 \int_{0}^{+\infty} \frac{ds}{2^{\lambda s/2} \<s\>^{3\alpha-1/4}}    \Big)\\
  &\leq &\dis   \vert a_0\vert+CT \Vert(a,g)\Vert^3_{T,\alpha}.
  \end{eqnarray*}
 Similarly, using that $\alpha \geq 1/4$
   \begin{eqnarray*}
 \sup_{s\geq 0} \<s\>^{\alpha}\vert \Phi_2(a,g)(t)\vert  &\leq& \sup_{0\leq t\leq T}\sup_{s\geq 0}( \<s\>^{\alpha}\vert  g_0(s)\vert)+CT \Big((\sup_{0\leq t\leq T} \vert a(t)\vert)^2\sup_{0\leq t\leq T}\sup_{s\geq 0}( \<s\>^{\alpha}\vert g(s))\vert     \\
 && +\big(\sup_{0\leq t\leq T} \vert a(t)\vert\big)\big(\sup_{0\leq t\leq T}\sup_{s\geq 0}( \<s\>^{\alpha}\vert g(s)\vert)\big)^2  + \big(\sup_{0\leq t\leq T}\sup_{s\geq 0} (\<s\>^{\alpha}\vert g(s)\vert)\big)^3     \Big)\\
  &\leq &\dis   \sup_{0\leq t\leq T}\sup_{s\geq 0}( \<s\>^{\alpha}\vert  g_0(s)\vert)+CT \Vert(a,g)\Vert^3_{T,\alpha}.
  \end{eqnarray*}
Thus, putting the previous estimates together and taking the sup norm in time, we obtain the bound
\begin{equation*}
 \Vert \Phi(a,g)\Vert_{T,\alpha} \leq  \Vert (a_0,g_0)\Vert_{T,\alpha} +CT \Vert (a,g)\Vert^3_{T,\alpha},
\end{equation*}
from which we deduce that $\Phi :X\longrightarrow X$ provided that $T\leq C  \Vert(a_0,g_0)\Vert^{-2}_{\alpha}$ with $C>0$ small enough. With similar estimates we can show that $\Phi$ is a contraction. 
\end{proof}

\subsection{Conserved quantities and symmetries}
\begin{itemize}
\item The "mass" $\displaystyle M = |a|^2 + \lambda \int_0^{+\infty} |g(s)|^2\,ds$.
\item The "kinetic energy" $\displaystyle E = \int_0^{+\infty} s |g(s)|^2\,ds$.
\end{itemize}
By Noether's theorem, this corresponds to the following symmetries, which leave the set of solutions invariant:
\begin{itemize}
\item Phase rotation: $(a,g) \mapsto (e^{i\theta} a \,,\, e^{i\theta} g)$, with $\theta \in \mathbb{R}$.
\item Phase modulation: $(a,g) \mapsto (a,e^{i \theta s} g)$, with $\theta \in \mathbb{R}$.
\end{itemize}
Another symmetry is given by the scaling
$$
(a(t),g(t,s)) \mapsto \lambda (a(\lambda^2 t),g(\lambda^2 t,s)), \qquad \mbox{for $\lambda > 0$}.
$$

\subsection{Transfer to high frequencies} The interactions between different frequencies, say $s$ and $2s$, are always damped by a coefficient $\frac{1}{2^{\lambda s}}$. It  suggests that any transfer of energy to high frequencies can only happen at a logarithmic speed - and this should also be the case for the full LLL equation.

\end{document}